\documentclass{article}

\newtheorem{theorem}{Theorem}
\newtheorem{lemma}{Lemma}
\newtheorem{corollary}{Corollary}

\newcommand{\bt}{\begin{theorem}}
\newcommand{\et}{\end{theorem}}
\newcommand{\bl}{\begin{lemma}}
\newcommand{\el}{\end{lemma}}
\newcommand{\bc}{\begin{corollary}}
\newcommand{\ec}{\end{corollary}}
\newcommand{\pf}{{\bf Proof}.\ }
\newcommand{\bq}{\begin{eqnarray*}}
\newcommand{\eq}{\end{eqnarray*}}
\newcommand{\be}{\begin{eqnarray}}
\newcommand{\ee}{\end{eqnarray}}
\newcommand{\beq}{\begin{equation}}
\newcommand{\eeq}{\end{equation}}
\newcommand{\benum}{\begin{enumerate}}
\newcommand{\eenum}{\end{enumerate}}
\newcommand{\ba}{\begin{array}}
\newcommand{\ea}{\end{array}}
\newcommand{\bfr}{\begin{flushright}}
\newcommand{\efr}{\end{flushright}}

\title{$N$--graphs, modular Sidon and sum--free sets,\\
and partition identities
}

\author{Melvyn B. Nathanson\thanks{
This work was supported in part by grants from the PSC--CUNY Research
Award Program and the NSA Mathematical Sciences Program.}\\
Department of Mathematics\\
Lehman College (CUNY)\\
Bronx, New York 10468\\
e-mail: nathansn@alpha.lehman.cuny.edu}
\date{}

\begin{document}
\maketitle

\begin{abstract}
Using a new graphical representation for partitions, the author obtains
a family of partition identities associated with
partitions into distinct parts of an arithmetic progression, or,
more generally, with partitions into distinct parts of a set that
is a finite union of arithmetic progressions associated
with a modular sum--free Sidon set.
Partition identities are also constructed for sets associated with
modular sum--free sets.
\end{abstract}

\section{$N$--graphs for partitions}

The {\em standard form} of a partition
$n = a_1 + a_2 + \cdots + a_k$ is
\[
\pi = (a_1,\ldots, a_k),
\]
where the parts $a_1,\ldots, a_k$ are positive integers arranged
in descending order.
The standard form of a partition is unique.

Associated to a partition $\pi = (a_1,\ldots, a_k)$ of $n$
is an array of dots, called the {\em Ferrers graph} of $\pi$.
This consists of $n$ dots arranged in $k$ rows,
with $a_1$ dots on the first row, $a_2$ dots on the second row,\ldots,
and $a_k$ dots on the $k$--th row.
The rows are aligned on the left.
The {\em Durfee square} $D(\pi)$ of the graph is the largest square array
of dots that appears in the upper left corner of the Ferrers graph.
We denote by $d(\pi)$ the number of dots on a side of the Durfee square,
or, equivalently, the number of dots on a diagonal of $D(\pi)$.

The Ferrers graph can be decomposed into a disjoint union of right angles,
called {\em hooks}.  The corners of the hooks are the dots on the
main diagonal of the Durfee square, and so the number of hooks is
$d(\pi)$.
Counting the number of dots on the hooks, we obtain the {\em hook numbers}
of the partition $\pi$.  Since the sum of the hook numbers is $n$,
we obtain a new partition of $n$, denoted $h(\pi)$ and called the
{\em hook number partition}.

For example, the partition $\pi = (7,6,6,5,4)$ of 28
has $d(\pi) = 4$.
The hook number partition is $h(\pi) = (11,8,6,3)$.

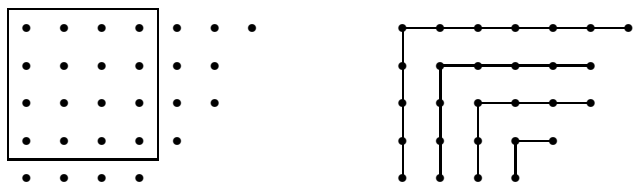
\begin{picture}(16,6)
\multiput(0,1)(1,0){4}{\circle*{0.2}}
\multiput(0,2)(1,0){5}{\circle*{0.2}}
\multiput(0,3)(1,0){6}{\circle*{0.2}}
\multiput(0,4)(1,0){6}{\circle*{0.2}}
\multiput(0,5)(1,0){7}{\circle*{0.2}}
\put(-0.5,5.5){\line(1,0){4}}
\put(-0.5,5.5){\line(0,-1){4}}
\put(3.5,5.5){\line(0,-1){4}}
\put(-0.5,1.5){\line(1,0){4}}

\multiput(10,1)(1,0){4}{\circle*{0.2}}
\multiput(10,2)(1,0){5}{\circle*{0.2}}
\multiput(10,3)(1,0){6}{\circle*{0.2}}
\multiput(10,4)(1,0){6}{\circle*{0.2}}
\multiput(10,5)(1,0){7}{\circle*{0.2}}
\put(10,5){\line(1,0){6}}
\put(10,5){\line(0,-1){4}}
\put(11,4){\line(1,0){4}}
\put(11,4){\line(0,-1){3}}
\put(12,3){\line(1,0){3}}
\put(12,3){\line(0,-1){2}}
\put(13,2){\line(1,0){1}}
\put(13,2){\line(0,-1){1}}
\end{picture}

MacMahon~\cite{macm23} introduced a beautiful arithmetic generalization
of the Ferrers graph of a partition.
Let $\pi = (a_1,\ldots, a_k)$ be a partition of $n$ in standard form.
For each positive integer $m$ we shall construct
the MacMahon {\em modular $m$--graph} of the partition $\pi$.
By the division algorithm, we can write each part $a_i$
uniquely in the form
\[
a_i  = u(a_i)m+s(a_i) \qquad\mbox{where $u(a_i) \geq 0$
and $1 \leq s(a_i) \leq m$.}
\]
The $m$--graph of $\pi$ consists of $k$ rows.
The $i$--th row has $u(a_i) + 1$ entries, where
the first $u(a_i)$ entries are $m$, and the last entry is $s(a_i)$.
In the special case $m=1$, we have $u(a_i) = a_i -1$ and $s(a_i) = 1$
for $i = 1,\ldots, k$.
The 1--graph is exactly the Ferrers graph with each dot replaced by 1.

The Durfee square of the $m$--graph, denoted $D_m(\pi)$,
is the largest square array of integers
contained in the upper left corner of the graph.
The number of dots a side of the Durfee square is denoted $d_m(\pi)$.
The hook number partition associated with the $m$--graph
is the partition $h_m(\pi)$ obtained by adding the numbers
on the hooks of the $m$--graph.
The hook number partition has $d_m(\pi)$ parts.

For example, if $\pi = (9,8,6,4)$, then
the $m$--graphs of $\pi$ for $m = 1, 2,$ and 3 are

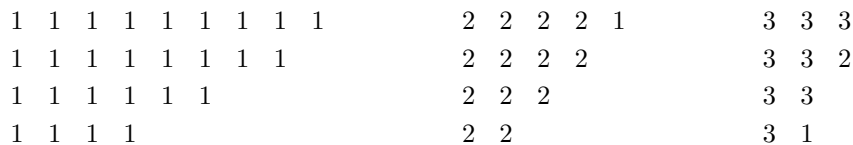
\begin{picture}(22,6)
\multiput(0,1)(1,0){4}{1}
\multiput(0,2)(1,0){6}{1}
\multiput(0,3)(1,0){8}{1}
\multiput(0,4)(1,0){9}{1}

\multiput(12,1)(1,0){2}{2}
\multiput(12,2)(1,0){3}{2}
\multiput(12,3)(1,0){4}{2}
\multiput(12,4)(1,0){4}{2}
\put(16,4){1}

\put(21,1){1}
\put(20,1){3}
\multiput(20,2)(1,0){2}{3}
\multiput(20,3)(1,0){2}{3}
\put(22,3){2}
\multiput(20,4)(1,0){3}{3}
\end{picture}

\noindent
Note that $d_1(\pi) = 4$, $d_2(\pi) = 3$, and $d_3(\pi) = 2$.
The hook number partitions associated with these graphs
are $h_1(\pi) = (12,9,5,1)$, $h_2(\pi) = (15,10,2)$,
and $h_3(\pi) = (18,9)$.

In this paper I introduce a generalization of MacMahon's $m$--graphs.
Let $m \geq 2$, and let $S = \{s_1,\ldots,s_{\ell}\}$ be
an ordered, nonempty set of positive integers
that are pairwise incongruent modulo $m$.
We do not assume that $1 \leq s \leq m$ for $s \in S.$
Let $A$ be the set of integers of the form
$um+s$, where $u \geq 0$ and $s \in S$.
Then $A$ is a finite union
of arithmetic progressions with difference $m$,
and every element $a \in A$ has a unique representation in the form
$a = u(a)m+s(a)$, where $u(a)$ is a nonnegative integer and $s(a) \in S$.

A partition $\pi$ into parts belonging to $A$
can be written uniquely in the form
\[
\pi = (a_1,\ldots, a_k)_N,
\]
where
\[
a_i = u(a_i)m+s(a_i) \in A,
\]
\[
u(a_1) \geq u(a_2) \geq \cdots \geq u(a_k),
\]
and, if $u(a_i) = u(a_{i+1}), s(a_i) = s_{j_i},$
and $s(a_{i+1}) = s_{j_{i+1}}$,
then $j_i \leq j_{i+1}$.
We shall call this the {\em standard $N$--form} for a partition
with parts in the set $A$.
Note that if $(a_1,\ldots, a_k)_N$ is the standard $N$--form
of a partition, then it is not necessarily true
that $a_i \geq a_{i+1}$ for all $i = 1,\ldots, k-1$.

The {\em $N$--graph} of the partition $\pi = (a_1,\ldots, a_k)_N$
will consist of $k$ rows.  The $i$--th row has $u(a_i) + 1$ entries,
where the first $u(a_i)$ entries are $m$ and the last entry is $s(a_i)$.
In particular, we obtain MacMahon's modular $m$--graphs
for partitions in the special case $\ell = m$ and
$s_j = m+1-j$ for $j = 1,\ldots, m$.

The Durfee square $D_N(\pi)$ of the $N$--graph
is the largest square array of integers
contained in the upper left corner of the graph,
and the hook number partition $h_N(\pi)$ associated with the $N$--graph
is the partition obtained by adding the numbers on the hooks
of the $N$--graph.  If $d_N(\pi)$ is the number of integers
on a side of the Durfee square $D_N(\pi)$,
then the hook number partition has $d_N(\pi)$ parts.

For example, let $m = 13$ and $S = \{3,2,20\},$ where
$\ell = 3$, and $s_1 = 3, s_2 = 2, s_3 = 20$.
Then
\bq
A & = & \{3, 16,29, 42,55,\ldots \} \cup \{2,15,28,41,54,\ldots \}
\cup \{20,33,46,59,\ldots\} \\
& = & \{2,3,15,16,20,28,29,33,41,42,46,54,55,59,\ldots\}.
\eq
Consider the partition $193 = 55 + 41 + 33 + 29 + 20 + 15.$
The standard form for this partition is
\[
\pi = (55,41,33,29,20,15),
\]
and the standard $N$--form is
\[
\pi = (55, 41, 29, 15, 33, 20)_N.
\]
The corresponding $N$--graph is

\setlength{\unitlength}{0.65cm}

\begin{picture}(6,8)
\multiput(0,6)(1,0){4}{13}
\put(4,6){3}
\multiput(0,5)(1,0){3}{13}
\put(3,5){2}
\multiput(0,4)(1,0){2}{13}
\put(2,4){3}
\put(0,3){13}
\put(1,3){2}
\put(0,2){13}
\put(1,2){20}
\put(0,1){20}
\put(-0.25,6.75){\line(1,0){3}}
\put(-0.25,6.75){\line(0,-1){3}}
\put(2.75,6.75){\line(0,-1){3}}
\put(-0.25,3.75){\line(1,0){3}}

\end{picture}

\noindent
The Durfee square contains 9 points, $d_N(\pi) = 3$,
and the hook number partition associated with the $N$--graph
is $h_N(\pi) = (127,63,3)$.

\section{Sum--free Sidon sets}
Let $S = \{s_1,\ldots, s_{\ell}\}$ be a nonempty finite set of integers,
and let $2S = \{s+s': s, s' \in S \}$.
The set $S$ is a {\em sum--free} if $S \cap 2S = \emptyset$.
The set $S$ is a {\em Sidon set} if every integer
has at most one representation as a sum of two elements of $S$,
that is, $s_{i_1} + s_{i_2} = s_{j_1} + s_{{j_2}}$
if and only if $\{{i_1},{i_2}\} = \{{j_1},{j_2}\}$.
For example, $\{1,6, 19\}$ is a sum--free Sidon set,
and $\{s\}$ is a sum--free Sidon set for every $s \neq 0$.

Let $m \geq 2$.
The set $S$ is {\em sum--free modulo $m$}
if the elements of $S$ are pairwise incongruent modulo $m$
and the congruence
$s_{i_1} + s_{i_2} \equiv s_{j}\pmod{m}$  has no solution with
$s_{i_1}, s_{i_2}, s_{j}\in S$.
The set $S$ is a {\em Sidon set modulo $m$} if every congruence class
modulo $m$ has at most one representation
as a sum of two elements of $S$,
that is, $S$ is a set of pairwise incongruent integers such that
$s_{i_1} + s_{i_2} \equiv s_{j_1} + s_{j_2} \pmod{m}$
if and only if $\{{i_1},{i_2}\} = \{{j_1}, {j_2}\}$.
For example, $\{1,6,19\}$ is a sum--free Sidon set modulo 15,
but not modulo 11, since $6 + 6 \equiv 0 + 1 \pmod{11}$.
The set $\{s\}$ is a sum--free Sidon set modulo $m$ for every
integer $s$ and every modulus $m$ that does not divide $s$.

Let $m \geq 2$, and let $S= \{s_1,\ldots, s_{\ell}\}$
be a set of positive integers that is a sum--free Sidon set modulo $m$.
Associated with $S$ is the set $A$ of positive integers of the form
$um+s$, where $u \geq 0$ and $s \in S.$
If $a = um+s \in A$, we define $u(a) = u$ and $s(a) = a$.
The integers $u(a)$ and $s(a)$ are uniquely determined by $a$.
For every positive integer $n$, let $\mathcal{A}(n)$
denote the set of partitions of $n$
in the form $n = a_1 + \cdots + a_k$, where $a_i \in A$ and
$u(a_i) > u(a_{i+1})$ for $i = 1,\ldots, k-1$.
Then $\pi = (a_1,\ldots,a_k)_N$ is the standard $N$--form of the partition.
Let $\mathcal{H}(n)$ denote the set of hook number partitions
associated with the $N$--graphs of the partitions in $\mathcal{A}(n)$.
The map that sends $\pi \in \mathcal{A}(n)$ to the hook number partition
$h_N(\pi) \in \mathcal{H}(n)$ is not, in general, one--to--one.
For example, let $m = 15$ and $S = \{1,6,19\}$.
The partitions $\pi^{(1)} = (96,61,64,21)_N$ and
$\pi^{(2)} = (96,66,64,16)_N$ have the same hook number partition
$h_N(\pi^{(1)}) = h_N(\pi^{(2)}) = (141,67,34).$
The standard $N$--graphs of the partitions $\pi^{(1)}$ and $\pi^{(2)}$
are

\begin{picture}(15,6)
\multiput(0,4)(1,0){6}{15}
\put(6,4){6}
\multiput(0,3)(1,0){4}{15}
\put(4,3){1}
\multiput(0,2)(1,0){3}{15}
\put(3,2){19}
\put(0,1){15}
\put(1,1){6}

\multiput(9,4)(1,0){6}{15}
\put(15,4){6}
\multiput(9,3)(1,0){4}{15}
\put(13,3){6}
\multiput(9,2)(1,0){3}{15}
\put(12,2){19}
\put(9,1){15}
\put(10,1){1}
\end{picture}

\noindent

Even though the map $\pi \mapsto h_N(\pi)$ is not one--to--one, there
is a partition identity that relates the sets $\mathcal{A}(n)$
and $\mathcal{H}(n)$.

\bt         \label{pari:theorem:Sidon}
Let $m \geq 2$ and let
\[
S= \{s_1,\ldots, s_{\ell}\}
\]
be a set of positive integers that is a sum--free Sidon set modulo $m$.
Let
\[
A = \{um+s :  u \geq 0 \mbox{ and } s \in S\}.
\]
Let $\mathcal{A}(n)$ be the set of partitions of $n$
in the form
\[
n = a_1 + \cdots + a_k,
\]
where
\beq         \label{pari:Sidon1}
a_i = u(a_i)m + s(a_i) \in A
\eeq
and
\beq         \label{pari:Sidon2}
u(a_1) > \cdots > u(a_k) \geq 0.
\eeq
Let $p_{\mathcal{A}}(n)$ denote the number of partitions
in the set $\mathcal{A}(n)$.

Let
\[
B = \{vm+s+s' : v \geq 1 \mbox{ and } s,s' \in S \}
\]
and
\[
H = A \cup B.
\]
Since $S$ is a sum--free Sidon set modulo $m$, each element $h \in H$
can be written uniquely in the form
\[
h = v(h)m + t(h),
\]
where $v(h) \geq 0$ and $t(h) \in S \cup 2S$.
Let $\mathcal{H}(n)$ be the set of partitions of $n$ of the form
\[
\pi' = (h_1,\ldots, h_d),
\]
where
\[
h_i \in H \qquad\mbox{for $i = 1,\ldots, d$},
\]
\[
v(h_i) - v(h_{i+1}) \geq 3 \qquad\mbox{for $i = 1,\ldots, d-1$,}
\]
and
\[
v(h_i) - v(h_{i+1}) \geq 4 \quad \mbox{if $h_{i+1} \in B$.}
\]
Let
\[
B' = \{vm+s+s' : v \geq 1 \mbox{ and } s,s' \in S, s \neq s' \}.
\]
For each partition $\pi' = (h_1,\ldots, h_d) \in \mathcal{H}(n)$,
let $e'(\pi')$ denote the number of $i \in \{1,\ldots, d\}$
such that $h_i \in B'$.
Then
\beq    \label{pari:Sidonidentity}
p_{\mathcal{A}}(n) = \sum_{\pi' \in \mathcal{H}(n)} 2^{e'(\pi')}.
\eeq
\et

\pf
Let $\pi = (a_1, \ldots, a_k)_N$ be the standard $N$--form
of a partition in $\mathcal{A}(n)$.
Then $a_1,\ldots, a_k$ are elements of the set $A$
that satisfy conditions~(\ref{pari:Sidon1}) and~(\ref{pari:Sidon2}).
Let $h_N(\pi)$ be the hook number partition
determined by the $N$--graph of $\pi$.
Then
\[
h_N(\pi) = (h_1,\ldots, h_d),
\]
where $d = d_N(\pi)$ is the number of integers on the side
of the Durfee square of the $N$--graph of $\pi$.
We shall show that $h_N(\pi)$ is a partition in $\mathcal{H}(n)$.

Each hook in the $N$--graph of $\pi$ consists of a horizontal row
of numbers and a vertical column of numbers;
the corner of the hook lies on the diagonal of the Durfee square.
The row consists of a sequence of $m$'s, and ends with an element of $S$.
The column consists of a sequence of $m$'s, and ends either with an $m$
or with an element of $S$.
In the first case the hook number is an element of $A$;
in the second case the hook must contain an $m$ on the diagonal,
and the hook number is an element of $B$.
Therefore, each hook number in the partition $h_N(\pi)$
belongs to the set $H = A \cup B.$

For $i = 1,\ldots, d$, we let $x_{i}$ denote the number of integers
on the row of the $i$--th hook,
and $y_{i}$ denote the number of integers
in the column below the corner of the $i$--th hook.
Let $1 \leq i \leq d -1$.
Since each row ends in an element of $S$, it follows
from~(\ref{pari:Sidon2}) that
\[
x_{i+1} \leq x_{i} - 2,
\]
and so the row in hook $i$ contains at least
two more elements equal to $m$ than the row in hook $i +1$.
Similarly,
\[
y_{i+1} \leq y_{i} -1,
\]
and the column in hook $i$ contains at least one
more element equal to $m$ than the column in hook $i + 1$.
Therefore, $v(h_{i}) - v(h_{i + 1}) \geq 3.$
If $h_{i +1} \in B$, then
the column of hook $i +1$ ends in an element of $S$,
hook $i$ has an $m$ to the left of this number,
and $v(h_{i}) - v(h_{i + 1}) \geq 4$.
This proves that the map $\pi \mapsto h_N(\pi)$ sends a partition
in $\mathcal{A}(n)$ to a partition in $\mathcal{H}(n)$.

Let $\pi' = (h_1,\ldots, h_d)$ be a partition in $\mathcal{H}(n)$
and let $e'(\pi')$ denote
the number of integers $i \in \{1, \ldots, d\}$ such that $h_i \in B'$.
We shall prove that there exist exactly $2^{e'(\pi')}$ partitions
$\pi \in \mathcal{A}(n)$ such that $h_N(\pi) = \pi',$
and we shall explicitly construct these partitions.

The partition $\pi' = (h_1,\ldots, h_d) \in \mathcal{H}(n)$ immediately
determines the shape of the $N$--graph of any partition
$\pi \in \mathcal{A}(n)$ such that $h_N(\pi) = \pi'$.
First, the Durfee square $D_N(\pi)$ must satisfy $d_N(\pi) = d$.
Second, we let $e$ denote the number of hook numbers
$h_i$ that belong to $B$.  Each of these hook numbers
is of the form $vm + s + s'$, where $s, s' \in S$, and the corresponding
hook in the $N$--graph of $\pi$ must contain two elements of $S$.
Each of the remaining $d-e$ hook numbers is of the form $vm + s$,
and the corresponding hook in the $N$--graph of $\pi$
contains only one element of $S$.
Therefore, the $N$--graph of $\pi$ contains
\[
2e + (d-e) = d+e = k
\]
elements of $S$.
Since each row of the $N$--graph of a partition contains
exactly one element of $S$, it follows that
the partition $\pi$ must contain exactly $k$ parts.
Thus, if $\pi' = h_N(\pi)$, then $\pi'$ determines
the number of parts in $\pi$.

Corresponding to the $e$ hook numbers $h_i \in B$ are
integers $1 \leq j_1 < \cdots < j_e \leq d$
such that $h_{j_i} \in B$ for $i = 1,\ldots, e$.
Then row $k$ in the standard $N$--graph of $\pi$
consists of $j_1 -1$ entries equal to $m$ followed by an element of $S$.
Similarly, row $k-1$ in the standard $N$--graph of $\pi$
consists of $j_2 -1$ entries equal to $m$ followed by an element of $S$.
In general, for $i = 0,1, \ldots, e-1$,
row $k-i$ in the standard $N$--graph of $\pi$
consists of $j_{i+1} -1$ entries equal to $m$ followed by an element of $S$.
This determines the shape of the bottom $e$ rows of the $N$--graph.
Then the hook numbers $h_1,\ldots, h_d$ determine the shape of the
top $d$ rows of the $N$--graph.
The only ambiguity concerns the elements of $S$ that are
at the ends of the rows.
If $h_i \equiv s\pmod{m}$ for some $s \in S$, then
the integer at the right end of row $i$ is $s$.
If $h_i \equiv 2s\pmod{m}$ for some $s \in S$, then
the integer at the right end of row $i$ is $s$ and
the integer at the bottom of column $i$ is $s$.
If $h_i \in B'$ and $h_i \equiv s+s'\pmod{m}$
for $s,s' \in S$ with $s \neq s'$,
then either the integer at the right end of row $i$ is $s$ and
the integer at the bottom of column $i$ is $s'$, or
the integer at the right end of row $i$ is $s'$ and
the integer at the bottom of column $i$ is $s$.
The set $S$ is a Sidon set modulo $m$, and so
these are the only ways to put elements of $S$
at the ends of the $i$--th hook
of the $N$--graph of $\pi$ to obtain the hook number $h_i$.
Since there are $e'(\pi')$ hook numbers $h_i$ that belong to $B'$,
it follows that there are exactly $2^{e'(\pi')}$
partitions $\pi \in \mathcal{A}(n)$
such that $h(\pi) = \pi'$.
This completes the proof.

\bt         \label{pari:theorem:AP}
Let $m \geq 2$ and let $s$ be a positive integer not divisible by $m$.
Let
\[
A = \{um+s :  u \geq 0 \},
\]
\[
B = \{vm+2s : v \geq 1 \},
\]
and
\[
H = A \cup B.
\]
Let $q_A(n)$ denote the number of partitions of $n$ as a sum of distinct
elements of $A$.
Let $q_H(n)$ denote the number of partitions of $n$ in the form
\[
n = h_1 + \cdots + h_d,
\]
where
\[
h_i = v(h_i)m + t(h_i) \in H \qquad\mbox{for $i = 1,\ldots, d$},
\]
\[
t(h_i) \in \{s,2s\},
\]
\[
v(h_i) - v(h_{i+1}) \geq 3 \qquad\mbox{for $i = 1,\ldots, d-1$,}
\]
and
\[
v(h_i) - v(h_{i+1}) \geq 4 \quad \mbox{if $h_{i+1} \in B$.}
\]
Then
\[
q_A(n) = p_H(n).
\]
\et

\pf
This follows immediately from Theorem~\ref{pari:theorem:Sidon},
applied to the sum--free Sidon set $S = \{s\}$ modulo $m$.

In the special case $m = 2$ and $s=1$, the set $A$ consists of all
odd positive numbers, and we obtain the following result
of Alladi~\cite{alla99}:
The number of partitions $n$ into distinct odd parts is equal
to the number of partitions of the form $n = h_1 + \cdots + h_d$,
where $h_d \neq 2$, $h_i - h_{i+1} \geq 6$ for $i = 1, \ldots, d-1$,
and $h_i - h_{i+1} \geq 7$ if $h_{i+1}$ is even.

\section{Partition identities for sum--free sets}
In the proof of Theorem~\ref{pari:theorem:Sidon},
the assumption that the sum--free set $S$ was a Sidon set
modulo $m$ implied that the cardinality
of the ``inverse image'' of a hook number
$h_i$ was at most two.
This produced the simple form
of the partition identity~(\ref{pari:Sidonidentity}).
We can also derive partition identities for sets $A$ that are
finite unions of arithmetic progressions constructed
from certain sets $S$ that are sum--free modulo $m$,
but not necessarily Sidon sets modulo $m$.
For example, we can consider sets $S$ that are sum--free modulo $m$
and have the property that if $s_{i_1},s_{i_2},s_{i_3},s_{i_4} \in S$ and
$s_{i_1} + s_{i_2} \equiv s_{i_3} + s_{i_4} \pmod{m}$,
then $s_{i_1} + s_{i_2} = s_{i_3} + s_{i_4}$.
The set of all odd  numbers $s$ such that $1 \leq s \leq m/2$
has this property.

\bt         \label{pari:theorem:sumfree}
Let $m \geq 2$ and let $S$
be a set of positive integers that is sum--free modulo $m$
and has the property that if $s_{i_1},s_{i_2},s_{i_3},s_{i_4} \in S$ and
$s_{i_1} + s_{i_2} \equiv s_{i_3} + s_{i_4} \pmod{m}$,
then $s_{i_1} + s_{i_2} = s_{i_3} + s_{i_4}$.
Let
\[
A = \{um+s :  u \geq 0 \mbox{ and } s \in S\}.
\]
Let $\mathcal{A}(n)$ be the set of partitions of $n$
in the form
\[
n = a_1 + \cdots + a_k,
\]
where
\[
a_i = u(a_i)m + s(a_i) \in A
\]
and
\[
u(a_1) > \cdots > u(a_k) \geq 0.
\]
Let $p_{\mathcal{A}}(n)$ denote the number of partitions in the set $\mathcal{A}(n)$.

Let
\[
B = \{um+s+s' : u \geq 1 \mbox{ and } s,s' \in S \}
\]
and
\[
H = A \cup B.
\]
Each element $h \in H$ can be written uniquely in the form
\[
h = v(h)m + t(h),
\]
where $v(h) \geq 0$ and $t(h) \in S \cup 2S$.
Let $\mathcal{H}(n)$ denote the set of partitions of $n$ of the form
\[
\pi' = (h_1,\ldots, h_d),
\]
where
\[
h_i  \in H \qquad\mbox{for $i = 1,\ldots, d$},
\]
\[
v(h_i) - v(h_{i+1}) \geq 3 \qquad\mbox{for $i = 1,\ldots, d-1$,}
\]
and
\[
v(h_i) - v(h_{i+1}) \geq 4 \quad \mbox{if $h_{i+1} \in B$.}
\]
For $h \in H$, let $r(h)$ denote the number of representations of $h$
as a sum of two elements of $S$, that is,
$r(h)$ is the number of ordered pairs $(s,s')$ such that $s+s'=h$ and
$s,s' \in S$.
Then
\[
p_{\mathcal{A}}(n)
= \sum_{\pi' \in \mathcal{H}(n) \atop \pi'= (h_1,\ldots, h_d)}
\prod_{i=1}^d r(h_i).
\]
\et

\pf
The proof is the same as the proof of Theorem~\ref{pari:theorem:Sidon}.

\end{document}